\documentclass[a4paper]{amsart}
\usepackage{amsmath, amssymb}
\usepackage{amssymb,pb-diagram}
\usepackage{amscd}

\usepackage{fancybox}
\usepackage{amsthm}

\usepackage{algorithm}
\usepackage{algorithmic}

\DeclareMathOperator{\divisor}{div}

\newcommand{\ZZ}{\mathbb Z}
\newcommand{\PP}{\mathbb P}

\newcommand{\CC}{\mathbb C}

\newcommand{\ded}{\mathfrak {d}}

\newcommand{\hor}{{\mathrm {hor}}}
\newcommand{\ver}{{\mathrm {ver}}}

\newcommand{\MW}{\mathop {\rm MW}\nolimits}

\newcommand{\Pic}{\mathop {\rm Pic}\nolimits}

\newcommand{\Supp}{\mathop {\rm Supp}\nolimits}

\newtheorem{thm}{Theorem}[section]
\newtheorem{mthm}{Theorem}
\newtheorem{cor}{Corollary}[section]
\newtheorem{prop}{Proposition}[section]
\newtheorem{lem}{Lemma}[section]
\newtheorem{defin}{Definition}[section]
\newtheorem{exmple}{Example}[section]
\newtheorem{rem}{Remark}[section]
\newtheorem{mrem}{Remark}
\newtheorem{qz}{Question}[section]

\newtheorem{mprbm}{Problem}

\newenvironment{example}{\begin{exmple}\rm }{\end{exmple}}
\newenvironment{remark}{\begin{rem}\rm }{\end{rem}}

\newcommand{\proofend}{\qed \par\smallskip\noindent}

\renewcommand{\thesubparagraph}{\theparagraph.\@arabic\c@subparagraph}
 
\begin{document}
  
\bigskip

\begin{center}

{\Large \bf  A remark on a Nagell-Lutz type statement for the Jacobian of a curve of genus $2$ 
 and  a 
  $(2, 3, 6)$ quasi-torus decomposition of a sextic with $9$ cusps}

\bigskip

Hiro-o TOKUNAGA and Yukihiro UCHIDA

\end{center}
{\bf Abstract.} In this note, an analogous statement to the Nagell-Lutz theorem does not hold for the Jacobian of a 
certain curve of genus $2$ over $\CC(t)$. As a by-product, we give a $(2, 3, 6)$ quasi-torus
decomposition for the dual curve of a smooth cubic.

\section*{Introduction}

Let $S$ be a smooth projective surface defined over $\CC$, the field of complex numbers.
 Let $\varphi: S \to C$ be an elliptic surface, i.e., $\varphi$ is relatively minimal and has a section $O$.
Under these circumstances, the generic fiber $E_S$ of $\varphi$ is an elliptic curve over $\CC(C)$, and
there exists a bijection between the set of $\CC(C)$-rational points $E_S(\CC(C))$ and the set of sections $\MW(S)$. Note
that  $O$ is regarded as the zero element in both group.
Suppose that $P \in E_S(\CC(C))$ is a torsion and $s_P$ denotes its corresponding section.  An analogous statement
to the Nagell-Lutz Theorem for an elliptic curve over $\ZZ$ can be stated:  \lq $s_P$ does not meet $O$.'
 It is known that this analogous statement is false when the characteristic of the base field is $p > 0$ and the order of $P$ is $p$
(see \cite[Appendix]{oguiso-shioda}).

 In this article we consider an analogous statement for a fibered surface whose general fiber is a curve of genus $2$.
 Let $\varphi : S \to \PP^1$ be such a surface. We call $\varphi$ {\it a genus $2$ fibration}, for short. 
 
 For a genus fibration $\varphi: S \to \PP^1$, we assume that $\varphi$ satisfies the following conditions throughout this article:
 
 \begin{itemize}
 
 \item  $\varphi$ is relatively minimal.
 
 \item Let $C_S$ be the generic fiber of $\varphi$.  $C_S$ is  a curve of genus $2$ over $\CC(t)$ given by
 an equation:
 \[
 y^2 = a_0x^6 + a_1x^5 + \ldots + a_6, \quad a_0 \in \CC^{\times}, \, a_i \,\, (i = 1, \ldots, 6) \in \CC[t] 
 \]
\item  $C_S$ has two rational points $\infty^{\pm}$, which give rise to two section $O^{\pm}$.

\end{itemize}

 Let $\Pic^0_{\CC(t)}(C_S)$ denotes the class group of $\CC(t)$-divisors of degree $0$.  It is known that
 any element of $\Pic_{\CC(t)}^0(C_S)$ can be uniquely presented of the form ${\mathfrak d} - (\infty^+ + \infty^-)$ such that
 ${\mathfrak d}$ is a $\CC(t)$-divisor of degree $2$ whose affine part is reduced in the sense of Definition~\ref{def:semi-reduced} (See \cite[Proposition 1]{ghmm}, for example).
  We then formulate our
 problem as follows:
 
 \begin{mprbm}\label{prbm:main}{Let  ${\mathfrak d} - (\infty^+ + \infty^-)$ be  a divisor of degree $0$ as above. Suppose that
 ${\mathfrak d} - (\infty^+ + \infty^-)$ gives rise to a torsion element in $\Pic_{\CC(t)}^0(C_S)$  
 Let  $D$ be a divisor on $S$ such that $D|_{C_S}$ gives the affine part of  ${\mathfrak d}$, i.e., points of ${\frak d} \neq \infty^{\pm}$.
  Put $D = D_{\hor} + D_{\ver}$, where
 $D_{\hor}$ is a divisor not containing any fiber component of $\varphi$ and $D_{\ver}$ consists of 
 fiber components. 
  Can we choose $D$ such that $\Supp(D_{\hor}) \cap (O^+ \cup O^-) = \emptyset$? 
 }
 \end{mprbm}

 This question asks if  a Nagell-Lutz type statement holds for $S$ or not. In this article, we first give a {\it negative} answer
 to Problem~\ref{prbm:main}.
  
 \begin{mthm}\label{thm:main}{There exists a genus $2$ fibration $\varphi : S \to \PP^1$ over $\PP^1$ such that
 
\begin{enumerate}
\item[(i)]  there exists a torsion element $\xi$ of order $3$ in $\Pic_{\CC(t)}^0(C_S)$,

\item[(ii)] let ${\mathfrak d}$ be the divisor of degree $2$ given in Proposition~\ref{prop:reduced-divisor1} so that
${\mathfrak d} - (\infty^+ + \infty^-)$ represents $\xi$, and

\item[(ii)] for  any divisor $D$ with $D|_{C_S} =$ the affine part of ${\mathfrak d}$, $\Supp(D_{\hor}) \cap (O^+ \cup O^-) \neq \emptyset$.

\end{enumerate}
 }
 \end{mthm}
 
 \begin{mrem}\label{mrem:main}{\rm If $C_S$ is given by
 \[
 y^2 = x^5 + \ldots + a_5, \quad a_i \in \CC[t]
 \]
 and ${\mathfrak d}$ is represented by a section $s$, then a Nagell-Lutz type statement holds by
 \cite{grant}.
 
 }
 \end{mrem}
 
 The surface $S$ in Theorem~\ref{thm:main} is realized as two point blowing-ups of the minimal resolution of
 a double cover of $\PP^2$ branched along a $9$-cuspidal sextic, i.e., a  sextic with $9$ cusps. A $9$-cuspidal
 sextic is the dual curve of a smooth cubic  and it has many $(2, 3)$ torus decompositions (see \cite{tokunaga99}, for 
 example). 
 Here a plane curve $B$ given by a homogeneous equation $F(T, X, Z) = 0$ is said to be a $(p,q)$  torus curve
(or a toric curve of type $(p,q)$), if $F$ can be represented as $F =  G^p + H^q$, where $G$ and $H$ are 
coprime homogeneous polynomials, respectively.  We call the right hand side a $(p, q)$ torus decomposition of 
$F$. Also $B$ is said to be
a  quasi-torus curve of type $(p, q, r)$ if $F$ satisfies $F_3^rF = F_1^p  + F_2^q$, where $F_1, F_2$ and $F_3$ are
pairwise coprime homogeneous polynomials. 
Torus curves have been studied by topological and arithmetic points of view by Cogolludo-Agustin, Kloosterman, Libgober, 
Oka and the first author (\cite{cogo-lib}, \cite{kloosterman2013, kloosterman2014}, \cite{oka}, \cite{oka-tai}, \cite{tokunaga99},
\cite{tokunaga99-2}).
 Among them, $(2,3)$ decompositions  for a $9$-cuspidal plane sextic $B$ is one of well-studied subjects. We here add one more observation
 to it: a $(2, 3, 6)$ quasi-torus decomposition of $B$ in Section~\ref{sec:torus}.

\medskip
\noindent
\textbf{Acknowledgements.}
The authors thank Professor Noboru Nakayama for his valuable comments.


\section{Preliminaries}

\subsection{A double cover construction for a genus $2$ fiber space over $\PP^1$}\label{subsec:double-cover}

We  fix some notation for later use. Let $\Sigma_d$ be the Hirzebruch surface of degree $d > 0$. We denote the 
section with self-intersection $-d$, the negative section, by $\Delta_0$ and a fiber of the ruling 
$\Sigma_d \to \PP^1$ by ${\mathfrak f}$.  A section linear equivalent to  $\Delta$ is denoted by $\Delta$. 
Likewise \cite{bannai-tokunaga15}, 
we take affine open subsets $U_1$ and $U_2$ of $\Sigma_d$ so that

\begin{enumerate}

 \item[(i)] $U_i \cong \CC^2$ $(i = 1, 2)$ and
 
 \item[(ii)] $U_1$ and $U_2$ have coordinates $(t, x)$ and $(s, x')$, respectively, with relations
 $t = 1/s, \, x = x'/s^d$.
 
 \end{enumerate}

Under these coordinates,  the negative section $\Delta_0$ is given by $x = x' = \infty$,  and 
a section $\Delta$ is given by $x = x' = 0$.

Let $f(t, x) \in \CC[t, x]$ be a polynomial of the form
\[
f(t, x):= a_0 x^6 + a_1(t)x^5 + \ldots + a_6(t), \quad a_i(t)   \in \CC[t],
\quad a_0 \in \CC^\times.
\]
If we choose a minimum positive integer $d$ such that $\deg a_i \le id$ ($i = 1, \ldots, 6$), then the affine curve $B^a$ given by
$f(t, x) = 0$ in $U_1$ gives a divisor $B$ on $\Sigma_d$ such that $B \sim 6\Delta, B\cap \Delta_0 = \emptyset$.

Choose $d \in \ZZ_{>0}$ as above. Let $f'_B : S'_B \to \Sigma_d$ be a double cover of $\Sigma_d$ with branch locus $B$ and let
 \[
\begin{CD}
 S'_{B} @<{\mu}<<S_{B} \\
           @VV{f'_{B}}V         @VV{f_B}V \\
\Sigma_d @<< {q}< \widehat{\Sigma_d}
\end{CD}
\]
be the diagram for the canonical resolution. The ruling on $\Sigma_d$ induces a fibration of curves of genus $2$, 
$\varphi_B : S_B \to \PP^1$, with sections $O^{\pm}$ which are the preimage of $\Delta_0$. The generic fiber
$C_B$ is given by
\[
y^2 = f(t, x).
\]
Hence $C_B$ is a curve of genus $2$ over $\CC(t)$. Note that $\mu : S_B \to S'_B$ is not the minimal resolution in general. 
For simplicity,   in the following, we always assume that

\medskip

\begin{center}
$(\ast)$ $\mu : S_B \to S'_B$ is the minimal resolution.
\end{center}

\medskip

By \cite[Lemma 5]{horikawa}, this condition is satisfied if $B$ has at worst simple singularities (see \cite{bpv} for simple singularities).

\subsection{A double cover branched along a $9$-cuspidal sextic}\label{subsec:double-sextic}

Let $B$ be a sextic with $9$ cusps. Let $f'_B :  S'_B \to \PP^2$ be a double cover of $\PP^2$ with branch locus $B$ and let
$\mu : S_B \to S'_B$ be the canonical resolution of $S'_B$ (see \cite{horikawa} for the canonical resolution) of $S'_B$ with
the following commutative diagram:
  \[
\begin{CD}
 S'_{B} @<{\mu}<<S_{B} \\
           @VV{f'_{B}}V         @VV{f_B}V \\
\PP^2 @<< {q}< \widehat{\PP^2},
\end{CD}
\]
where $q$ is a composition of $9$ time blowing-ups at the $9$ cusps. Choose a point $z_o \in \PP^2\setminus B$. By abuse of notation, we denote $q^{-1}(z_o)$ by
$z_o$ as $q$ is identity over $z_o$. Let  $q_{z_o} : (\widehat{\PP^2})_{z_o} \to \PP^2$ be a blowing-up at $z_o$ and 
let $S_{B, z_o}$ be the induced double cover of $(\widehat{\PP^2})_{z_o}$:
 \[
\begin{CD}
 S_{B} @<{\mu}<<S_{B, z_o} \\
           @VV{f_{B}}V         @VV{f_{B, z_o}}V \\
\widehat{\PP^2} @<< {q_{z_o}}< (\widehat{\PP^2})_{z_o}.
\end{CD}
\]
A pencil of lines through $z_o$ induces a fiber space of curve of genus $2$ on $S_{B, z_o}$, which we denote by
$\varphi_{B, z_o} : S_{B, z_o} \to \PP^1$.  We choose a homogeneous coordinates $[T, X, Z]$ of $\PP^2$ such that
$z_o=[0,1,0]$ and $B$ is given by
\[
F_B(T, X, Z) = a_0X^6 + a_1(T, Z)X^5 + \ldots + a_6(T, Z) = 0, a_0 \in \CC^{\times}
\]
where $a_i$ ($i = 1, \ldots, 6$) are homogeneous polynomials of degree $i$. Put $f_B(t,x) := F_B(t, x, 1)$. Under this 
circumstance, the generic fiber $C_{B, z_o}:= C_{S_{B, z_o}}$ of $\varphi_{B, z_o}$ is given by
\[
y^2 = f_B(t, x).
\]

Since the curve given by $f_B(t, x) = 0$  also defines a divisor $B_1$ on in $\Sigma_1$, we have a double cover $S_{B_1}$
of $\Sigma_1$. 
branched over $B_1$ as in Section~\ref{subsec:double-cover}.
Not that $S_{B, z_o} = S_{B_1}$ and $(q_{z_o}\circ f_{B,z_o})^{-1}(z_o) = O^+ + O^-$.

We show that there exists a $3$-torsion in $\Pic^0_{\CC(t)}(C_{B, z_o})$ described in Theorem~\ref{thm:main}.

\subsection{Divisors representing elements in $\Pic^0_K(C)$ for a hyperelliptic curve $C$}

Let $K$ be a field with ${\rm {char}}(K)=0$ and $\overline{K}$ denotes its algebraic closure.
As for this section, main references are \cite{ghmm, mwz, silverman}. In particular, we use the notation in 
\cite[Chapter II]{silverman}.

Let $C$ be a hyperelliptic curve of genus $g$ given by the equation
\[
y^2 = f(x)
:= a_0 x^{2g+2} + a_1x^{2g+1} + \ldots + a_{2g+2}, \quad a_i \in K,
\quad \sqrt{a_0} \in K^\times.
\]

 We denote the hyperelliptic involution of $C$ by $\sigma$.
 
 \begin{defin}\label{def:div-support}{\rm For a divisor ${\mathfrak d} = \sum_im_iP_i$, the support of ${\mathfrak d}$ is the set 
 $\Supp(D) :=\{ P_i \mid m_i \neq 0\}$.
 }
 \end{defin}
 
 \begin{defin}\label{def:semi-reduced}{\rm An effective divisor ${\mathfrak d} = \sum m_iP_i$ on $C$  is said to be semi-reduced
 if $P_i \neq \infty^{\pm}$ and $\sigma(P_i) \not\in \Supp({\mathfrak d})$, unless $P_i = \sigma(P_i)$ and $m_i =1$. Moreover,
 if $\sum_im_i \le g$, ${\mathfrak d}$ is said to be a reduced divisor.
 }
 \end{defin}
 
  \begin{rem}\label{rem:reduced}{\rm Note that a (semi-) reduced divisor in Definition~\ref{def:semi-reduced} is
  different  from a {\it reduced divisor} in usual sense, i.e., a divisor without multiple components.}
  \end{rem}
 
 \begin{prop}\label{prop:reduced-divisor1}{(cf. \cite[Propositon 1]{ghmm})Suppose that $g$ is even.  For a $K$-rational divisor ${\mathfrak d} \in \Pic^0_K(C)$, 
 ${\mathfrak d}$ has unique representative in $\Pic_K^0(C)$ of the form ${\mathfrak d}_0  - g/2(\infty^+ + \infty^-)$, where
 ${\mathfrak d}_0$ is an effective $K$-rational divisor of degree $g$ whose affine part is reduced.
 }
 \end{prop}

 Since a curve of genus $2$ is hyperelliptic, we have:
 
 \begin{cor}\label{cor:deduced-genus2}{If $g = 2$,
 for a $K$-rational divisor ${\mathfrak d} \in \Pic^0_K(C)$, 
 ${\mathfrak d}$ has unique representative in $\Pic_K^0(C)$ of the form ${\mathfrak d}_0  - (\infty^+ + \infty^-)$, where
 ${\mathfrak d}_0$ is an effective $K$-rational divisor of degree $2$ whose affine part is reduced, i.e.,
 ${\mathfrak d}_0$ is of the form $P_1 + P_2, P + \infty^{\pm}, 2\infty^+$ or $2\infty^-$, where $P_i, P$ are finite points and
 $P_2 \neq \sigma(P_1)$.
 }
 \end{cor}

\subsection{An algorithm for a reduction of a divisor}
 
We here recall an algorithm by which we can  compute
a reduced divisor from a given semi-reduced divisor  of higher degree after the one given in~\cite{ghmm}.

Let $C$ be a hyperelliptic curve of genus $2$ as above.
It is known that an effective affine semi-reduced divisor ${\frak d}$ on $C$
is represented as a pair of two polynomials,
which is called the Mumford representation of ${\frak d}$.
Mumford representations derive from the construction of
the Jacobian variety of a hyperelliptic curve given by Mumford~\cite{mumford}.
Let $\ded=\sum_i m_i P_i$ be an effective affine semi-reduced divisor.
Let $P_i=(x_i,y_i)$ and $u(x)=\prod_i (x-x_i)^{m_i}$.
Then there exists a unique polynomial $v(x)$ such that
$\deg(v)<\deg(u)$, $v(x_i)=y_i$, and $u\mid (v^2-f)$
(see the first half of~\cite[Theorem~A.5.1]{mwz},
where this fact is proved when $\deg(f)$ is odd.
The same proof, however,   works when $\deg(f)$ is even).
We call the pair $(u,v)$ the Mumford representation of $\ded$.
Conversely, let $(u,v)$ be a pair of polynomials such that
$u$ is monic and $v$ satisfies the above conditions.
Then there exists a unique effective affine semi-reduced
divisor $\ded$ such that
the Mumford representation of $\ded$ is equal to $(u,v)$.
We denote the divisor $\ded$ by $\divisor[u,v]$.

\begin{remark}
The definition of $\divisor[u,v]$ is different from
that of $\divisor(u,v)$ in~\cite{mwz}.
In~\cite{mwz}, $\deg(f)$ is assumed to be odd
and $\divisor(u,v)$ is defined by
$\divisor(u,v)=\ded - \deg(\ded) \infty$,
where $\infty$ is the point at infinity.
When $\deg(f)$ is even, $C$ has two points at infinity
and the behavior of the polar divisor of $y-v(x)$
is complicated as described later.
Thus we defined $\divisor[u,v]$ as an effective affine divisor
as in~\cite{ghmm}.
\end{remark}

From now on, we assume that $\deg(f)=6$.
Let $\ded_\infty = \infty^+ + \infty^-$.
For a given divisor class $\ded\in\Pic_K^0(C)$,
by Corollary~\ref{cor:deduced-genus2},
$\ded$ has a unique representative of the form 
$\ded_0-\ded_\infty$, where $\ded_0$ is an effective $K$-rational divisor
of degree $2$ whose affine part is reduced.
We can compute the divisor $\ded_0$ by the reduction algorithm,
which is a part of Cantor's algorithm~\cite{cantor}.
Although Cantor~\cite{cantor} assumed that $\deg(f)$ is odd,
we can generalize his algorithm to the case $\deg(f)$ is even.
We follow the description in~\cite{ghmm}.

Let $x'=1/x$, $y'=y/x^3$, and $\tilde{f}(x')=(x')^6 f(1/x')$.
Then $\tilde{f}(x')$ is a polynomial of degree $5$ or $6$
and we have $(y')^2=\tilde{f}(x')$.
Since $\deg(f)=6$, we have $\tilde{f}(0)\neq 0$.
We define
\[
	a_+ = y'(\infty^+) = \left(\frac{y}{x^3}\right)(\infty^+), \quad
	a_- = y'(\infty^-) = \left(\frac{y}{x^3}\right)(\infty^-).
\]
Then $a_+$ and $a_-$ are the square roots of $\tilde{f}(0)$,
hence we have $a_+\neq a_-$.
For a polynomial $p(x)$, 
we denote by $\operatorname{lt}(p)$ the leading term of $p(x)$.
If $\operatorname{lt}(p)=a_+ x^3$, then 
the function $y-p(x)$ has a pole of order less than $3$ at $\infty^+$
and a pole of order $3$ at $\infty^-$.
The case $\operatorname{lt}(p)=a_- x^3$ is similar.
If $\operatorname{lt}(p)\neq a_+ x^3, a_- x^3$, then
$y-p(x)$ has poles of order $\max\{3, \deg(p)\}$
at $\infty^+$ and $\infty^-$.

In the following, we only need to consider divisors $\divisor[u,v]$
of degree $4$ such that $\operatorname{lt}(v)\neq a_+ x^3, a_- x^3$.
For such divisors, the reduction algorithm is
described as in Algorithm~\ref{alg:reduction},
which is a specialization of \cite[Algorithm~2]{ghmm}.
In general, we can compute the reduced divisor associated with
a semi-reduced divisor by \cite[Algorithms~2 and~3]{ghmm}.

\begin{algorithm}
\caption{Reduction}
\label{alg:reduction}
\begin{algorithmic}[1]
\REQUIRE A semi-reduced affine divisor $\ded_0=\divisor[u_0,v_0]$ such that
$\deg(\ded_0)=4$ and $\operatorname{lt}(v_0)\neq a_+ x^3, a_- x^3$.
\ENSURE The reduced affine divisor $\ded_1=\divisor[u_1,v_1]$ such that
$\deg(\ded_1)=2$ and
$\ded_0-\ded_1\sim \ded_\infty$.
\STATE $u_1':=(f-v_0^2)/u_0$.
\STATE $u_1:=u_1'/\operatorname{lc}(u_1')$,
where $\operatorname{lc}(u_1')$ is the leading coefficient of $u_1'$.
\STATE $v_1:=-v_0\bmod u_1$.
\RETURN $\divisor[u_1,v_1]$.
\end{algorithmic}
\end{algorithm}

The geometric interpretation of the reduction algorithm
is also given in \cite{ghmm}.
We explain it for Algorithm~\ref{alg:reduction}.
Let $\ded_0=\divisor[u_0,v_0]$ as in Algorithm~\ref{alg:reduction}.
Then the zero divisor $\ded_z$ of the function $y-v_0(x)$ satisfies
$\ded_z=\ded_0+\sigma(\ded_1)$, where $\sigma(\ded_1)$ is the divisor
obtained by applying the hyperelliptic involution $\sigma$
to each point appearing in $\ded_1$.
Since $\deg(\ded_0)=4$ and $\operatorname{lt}(v_0)\neq a_+ x^3, a_- x^3$,
we have $\deg(\ded_z)=6$, hence $\deg(\ded_1)=2$.
Therefore we have
\[
	\left(\frac{y-v_0(x)}{u_0(x)}\right)
	= \sigma(\ded_1) - \sigma(\ded_0)
	+ \ded_\infty.
\]
Since $\ded+\sigma(\ded)\sim\deg(\ded)\ded_\infty$
for any divisor $\ded$, we have
$\ded_0-\ded_1\sim\ded_\infty$.


\section{Proof of Theorem~\ref{thm:main}}\label{sec:proof}

\subsection{An $S_3$-cover of $\PP^2$ branched along a $9$-cuspidal sextic}\label{subsec:S3-cover}

Let $E$ be an smooth cubic in $\PP^2$. We choose a flex $O$ of $E$ and fix it. Let $S^3(E)$ be the $3$-fold symmetric product, 
which is the set of effective divisors of degree $3$. Let $\alpha : E\times E\times E \to \Pic^0(E)$ be the
Abel-Jacobi map given by
 $(P_1, P_2, P_3) \mapsto P_1 + P_2 + P_3 - 3O$. Note that, as $E$ is identified with $\Pic^0(E)$ by $P \mapsto 
P - O$, $\alpha((P_1, P_2, P_3)) = P_1\dot{+}P_2\dot{+}P_3$, where $\dot{+}$ denotes the addition on $E$ with $O$ as the zero.
Let $\varpi: E\times E\times E \to S^3(E)$ be an $S_3$-cover given by the definition of $S^3(E)$ and let $\overline {\alpha} :
S^3(E) \to \Pic^0(E)$ be the induced map. By its definition, we have the following:

\begin{enumerate}

 \item[(i)] $A:=\alpha^{-1}(0)$ is an Abelian surface.
 
 \item[(ii)] $\overline{\alpha}^{-1}(0)$ is $\PP^2$.
 
 \item[(iii)] $\pi := \varpi|_A: A \to \PP^2$ is an $S_3$-cover of $\PP^2$ and 
 $\alpha = \overline{\alpha}\circ\pi$.

\end{enumerate}

Since $\overline{\alpha}^{-1}(0)$ is the set of effective divisors of degree $3$  cut out by lines, $\overline{\alpha}^{-1}(0)$ is considered
as the set of lines. Hence, we infer that the branch locus of $\pi$ is the dual curve of $E$, which is a $9$-cuspidal sextic.

\subsection{Proof}\label{subsec:proof}

Let $\pi : A \to \PP^2$ be the $S_3$-cover given in Section~\ref{subsec:S3-cover} and we put $\Delta_{\pi} = B$.
Let  $\beta_1(\pi): D(A/\PP^2) \to \PP^2$  and $\beta_2(\pi) : A \to D(A/\PP^2)$ be the  double and cyclic triple covers given
in \cite{tokunaga94}. (Note that $S_3 \cong D_6$).
$\beta_1(\pi)$ coincides with the double cover $f'_B : S'_B \to \PP^2$
considered in Section~\ref{subsec:double-sextic}, and $\beta_2(\pi): A \to S'_B$  is a cyclic triple cover branched over
$9$ singularities which are of type $A_2$ only. Hence it induces a cyclic triple cover $g_B : X_B \to S_B$ with
branch locus $\sum_{i=1}^9 (\Theta_{i,1} + \Theta_{i,2})$.

Hence by relabelling $\Theta_{i, j}$ $(i =1, \ldots, 9, j = 1, 2)$ suitably, by \cite[Lemma 8.7]{tokunaga02},
we may assume that there exists a divisor $D_o$ on $S_{B}$ such that
\[
(\ast) \quad 3D_o \sim \sum_{i=1}^9 (2\Theta_{i,1} + \Theta_{i,2}). 
\]
Choose $z_o \in \PP^2$ such that any line through $z_o$ passes through at most one cusp of $B$. Let  $\varphi_{B, z_o} : S_{B, z_o}
\to \PP^1$ be the genus $2$ fibration as in Section~\ref{subsec:double-sextic} and let $C_{B, z_o}$ be its generic fiber. 

\begin{lem}\label{lem:3-tor}{Put ${\mathfrak d}_o = D_o|_{C_{B, z_o}}$. Then the class $[{\mathfrak d}_o]$ in 
$\Pic_{\CC(t)}^0(C_{B,z_o})$ is a $3$-torsion.
}
\end{lem}

\proof Choose $\phi \in \CC(S_{B, z_o})$ such that 
\[
(\phi) = 3D_o - \sum_{i=1}^9(2\Theta_{i, 1} + \Theta_{i, 2}).
\]
This relation implies $[3{\mathfrak d}_o]$ is $0$ in $\Pic_{\CC(t)}(C_{B, z_o})$. We show that ${\mathfrak d}_o \not\sim 0$ in $\Pic_{\CC(t)}^0(C_{B, z_o})$.
If ${\mathfrak d}_o \sim 0$ on $C_{B, z_o}$, then ${\mathfrak d}_o$ is linear equivalent to $(h)$, $h \in \CC(t)^{\times}$ and 
we may assume that $D_o$ is linear equivalent to a divisor whose irreducible 
components are all in fibers of $\varphi_{B, z_o}$.
Put
\[
D_o \sim aF + \sum_{i=1}^9(b_{i, 1}\Theta_{i, 1} + b_{i, 2}\Theta_{i, 2}), a, b_{i, j} \in \ZZ.
\]
Then by $(\ast)$ we have
\[
3aF + \sum_{i=1}^9\{(3b_{i, 1} - 2)\Theta_{i, 1} + (3b_{v,2} -1)\Theta_{i, 2}\} \sim  0.
\]
By computing the intersection product  with $O^+$, we have
\[
0 = 3aF\cdot O^+, \, \mbox{i.e.,}\,  a = 0.
\]
and 
\[
\sum_{i=1}^9\{(3b_{i, 1} - 2)\Theta_{i, 1} + (3b_{v,2} -1)\Theta_{i, 2}\} \sim  0.
\]
We next compute the intersections with $\Theta_{i, 1}, \, \Theta_{i, 2}$ and we have
\begin{eqnarray*}
-2(3b_{i, 1} - 2) + (3b_{i, 2} - 1) & = & 0 \\
(3b_{i,1} -2)  -2 (3b_{i, 2} -1) & = & 0.
\end{eqnarray*}
This implies $b_{i, 1} = 2/3, b_{i,2} = 1/3$, which is impossible as $b_{i, j} \in \ZZ$.
\proofend

By Corollary~\ref{cor:deduced-genus2}, there exist an effective divisor ${\mathfrak d}$ such that $\deg {\mathfrak d} = 2$
and ${\mathfrak d} - {\mathfrak d}_{\infty} \sim {\mathfrak d}_o$.  Let $D$ be any divisor such that  $D|_{C_{B, z_o} }  = 
{\mathfrak d}$. Then as ${\mathfrak d}_o \sim {\mathfrak d} - {\mathfrak d}_{\infty}$, we have  
$D_o \sim D - O^+ - O^- + {\mathcal F}_o$, ${\mathcal F}_o$ is a vertical divisor. Put
 \[
D_{\ver} + {\mathcal F}_o = aF + \sum_{i=1}^9(b_{i, 1}\Theta_{i,1} + b_{i, 2}\Theta_{i,2}), \,\, a, b_{i,j} \in \ZZ
\]

\begin{lem}\label{lem:main}{$(\Supp(D_{\hor})\setminus \{O^{\pm}\})\cap (O^+\cup O^-) \neq \emptyset$.}
\end{lem}
\proof By Corollary~\ref{cor:deduced-genus2}, 
then  $D$ can be written one of the following forms
\[
D_{\hor} = \left \{ \begin{array}{c}
                    C \\
                    s + O^+ \\
                    s + O^- \\
                    2O^+ \\
                    2O^- 
                    \end{array} \right. ,
\]                    
where $C$ is a horizontal divisor such that (a) $C\cdot F = 2$ and (b) it does not contain $O^{\pm}$ as its irreducible components,  and   $s$ is a
 section $\neq O^{\pm}$.

\medskip

\underline{The case $D_{\hor} = C$}. By our assumption, $\overline{C}:= q\circ q_{z_o}\circ f_{B, z_o}(C)$ does not contain 
$z_o$. Hence $\overline{C}$  meets 
a line through $z_o$ at $2$ points and it implies $\overline {C}$ is a curve of degree $2$.

 By the relation $(\ast)$, we have
\[
3(C - O^+ - O^- + aF + \sum_{i=1}^9(b_{i, 1}\Theta_{i,1} + b_{i, 2}\Theta_{i,2}) ) \sim  \sum_{i=1}^9(2\Theta_{i, 1} + \Theta_{i, 2}).
\]

{\bf Claim.} $\displaystyle{\left [\begin{array}{c}
        C\cdot \Theta_{i,1} \\
       C\cdot \Theta_{i,2}  
       \end{array} \right ] \neq  \left [\begin{array}{c}
       0 \\
       0
       \end{array} \right ]}$.

       Proof of Claim. Suppose that  $\displaystyle{\left [\begin{array}{c}
        C\cdot \Theta_{i,1} \\
       C\cdot \Theta_{i,2}  
       \end{array} \right ] =  \left [\begin{array}{c}
       0 \\
       0
       \end{array} \right ]}$.  Then the above equivalence implies
       \begin{eqnarray*}
      && \left [ \begin{array}{c}
       3(C - O^+ - O^- + aF + \sum_{i=1}^9(b_{i, 1}\Theta_{i,1} + b_{i, 2}\Theta_{i,2}))\cdot \Theta_{i,1} \\
        3(C - O^+ - O^- + aF + \sum_{i=1}^9(b_{i, 1}\Theta_{i,1} + b_{i, 2}\Theta_{i,2}))\cdot \Theta_{i,2}
       \end{array} \right ] \\
       & =& 
        \left [ 
         \begin{array} {c}
         (\sum_{i=1}^9(2\Theta_{i, 1} + \Theta_{i, 2}))\cdot \Theta_{i,1} \\
       (\sum_{i=1}^9(2\Theta_{i, 1} + \Theta_{i, 2}))\cdot \Theta_{i,2} 
       \end{array}
        \right ] 
       \end{eqnarray*}
      Since $F \cdot\Theta_{i, j} = O^{\pm}\cdot\Theta_{i, j} = 0$, the left hand side is
      \[
      \left [ \begin{array}{c}
       -6b_{i,1} + 3b_{i,2} \\
       3b_{i,1} - 6b_{i,2}
       \end{array}
       \right ]
       \]
       while the right hand side is  $\displaystyle{\left [\begin{array}{c}
                           -3 \\
                           0
                           \end{array}
                       \right ]}$.  Hence $b_{i,1} = 2/3, \, b_{i,2} = 1/3$, which contradicts to $b_{i,j} \in \ZZ$. This shows Claim.

       \medskip
       
       By Claim,  we infer that $\overline{C}$ passes through the $9$ cusps of $B$. Since $B$ is a sextic, this is impossible.

       \medskip

      \underline{The case $D_{\hor} = s + O^+, s+ O^-$}.   In this case, $\overline{s}:= q\circ q_{z_o}\circ f_{B, z_o}(s)$ 
      does not contain $z_o$. Hence $\overline{s}$ intersects  a line through $z_o$ at one point and we infer that $\overline{s}$
      is a line. 
%
If  $\displaystyle{\left [\begin{array}{c}
        s\cdot \Theta_{i,1} \\
       s\cdot \Theta_{i,2}  
       \end{array} \right ] =  \left [\begin{array}{c}
       0 \\
       0
       \end{array} \right ]}$,  as  $F \cdot\Theta_{i, j} = O^{\pm}\cdot\Theta_{i, j} = 0$, a similar argument to the previous case shows that $b_{i, j} \not\in  \ZZ$. Hence 
$\displaystyle{\left [\begin{array}{c}
        s\cdot \Theta_{i,1} \\
       s\cdot \Theta_{i,2}  
       \end{array} \right ] \neq \left [\begin{array}{c}
       0 \\
       0
       \end{array} \right ]}$.   This implies that $\overline{s}$ passes through all $9$ cusps and $\overline{s}$ intersects $B$ at 
       $9$ points, which is impossible.    
       
       \medskip
       
        \underline{The case $D_{\hor} =  2O^+,  2O^-$}.  As $O^{\pm}\Theta_{i, j} = 0$, a similar argument to previous two cases
         shows that $b_{i, j} \not\in \ZZ$. Hence this case does not occur.
         
        By Lemmas~\ref{lem:3-tor} and~\ref{lem:main}, we have Theorem~\ref{thm:main}.
  

\section{A quasi-torus decompositions of type $(2,3,6)$ for  a $9$ cuspidal sextic}\label{sec:torus}
 
 In this section, we consider an explicit description of $D_o$ in Section~\ref{subsec:proof}. To this purpose, 
 we consider an explicit description
of $\pi : \to \PP^2$.  As an application, 
we give a $(2, 3, 6)$ quasi-torus decomposition  for $B$ explicitly,  Suppose that  $B$ as  a dual curve of an elliptic curve $E$ given by
an equation
\[
E: v^2 = u^3 + au + b, \quad 4a^3 + 27b^2 \neq 0.
\]
Let $l_{t,x}$ be a line given by $v = xu + t$. Then $(t, x) \in B$ if and only if $l_{t, x}$ is tangent to $E$. By equating $E$ and $l_{t,x}$ with respect to 
$v$, we have an cubic equation with respect to $u$:
\[
p(u):= {u}^{3}-{x}^{2}{u}^{2}- \left( 2\,tx-a \right) u-{t}^{2}+b = 0
 \]
The discriminant of $D_E(t, x)$ with respect to $u$ is given by
\begin{eqnarray*}
D_E(t, x) & = & \,-4\,at{x}^{5}+4\,b{x}^{6}+{a}^{2}{x}^{4}-4\,{t}^{3}{x}^{3}-30\,a{t}^{2}{x}^{2}+36\,bt{x}^{3} \\
 & & +24\,{a}^{2}tx-18\,ab{x}^{2}-27\,{t}^{4}- 4\,{a}^{3}+54\,b{t}^{2}
\mbox{}-27\,{b}^{2}.
\end{eqnarray*}

The dual curve $B := E^{\vee}$ is given by $F_B(T, X, Z) := Z^6D_E(T/Z, X/Z) = 0$ (see \cite[Ch. 5]{fischer}). Put
\begin{eqnarray*}
g &:= & -\frac 13x^4-2tx+a \\
h & := &  -{\frac {2}{27}}{x}^{6}- \frac 13 \left( 2\,tx-a \right) {x}^{2}-{t}^{2}+b.
\end{eqnarray*}
Then we have 
\[
D_E(t, x) := -4g^3 - 27 h^2.
\]
Hence we have a $(2, 3,6)$ quasi-torus decomposition  of $B$:
\begin{eqnarray*}
Z^6F_B(T, X, Z) & =& \left \{\sqrt[3]{-4}\left (-\frac13 X^4 - 2TXZ^2 + aZ^4\right)\right \}^3 \\
                  & &     + \left\{\sqrt{-27}\left (-{\frac {2}{27}}{X}^{6}- \frac 13 \left( 2\,TXZ^2-aZ^4 \right) {X}^{2}-{T}^{2}Z^4+bZ^6\right )\right \}^2.
\end{eqnarray*}       

We show that the $(2,3,6)$ quasi-torus decomposition gives rise to a $3$-torsion of $\Pic^0_{\CC(t)}(C_{B, z_o})$ considered in Section~\ref{sec:proof}.   

By our observation in Section~\ref{subsec:S3-cover}, we infer  that the rational function field $\CC(A)$ coincides with
the minimal splitting field  for the cubic equation $p(u) = 0$.  
By taking the Cardano formula on the cubic equation into account, we infer that
\begin{itemize}
 
 \item $\CC(D(A/\PP^2)) = \CC(t, x, \sqrt{D_E(t, x)})$, where $\CC(D(A/\PP^2))$
 is the rational function field of $D(A/\PP^2)$
 and 
 
 \item $\CC(A)$ is given by $\CC(D(A/\PP^2))\left (\sqrt[3]{h + \sqrt{D_E(t, x)}}\right )$.
 
 \end{itemize}
 
 Since $(h + \sqrt{D_E(t, x)})(h - \sqrt{D_E(t,x)}) = g^3$, from above fact, we infer that
 the rational function $h + \sqrt{D_E(t, x)}$ can be considered the rational function $\phi$ in our proof of 
 Lemma~\ref{lem:3-tor} and it gives a divisor on $S_B$ of the form
 \[
 (h + \sqrt{D_E(t, x)} ) = 3D_o - \sum_{i=1}^9(2\Theta_{i,1} + \Theta_{i, 2}), \quad D_o \not\sim 0.
 \]
 Hence, the class $D_o$ gives rise to a $3$-torsion in $\Pic_{\CC(t)}^0(C_{B, z_o})$ considered in Section~\ref{sec:proof}.
 In the next section, we consider  a plane curve  which gives rise to $D_{\hor}$  considered in
 Section~\ref{sec:proof}.

 \section{Examples}\label{sec:examples}
In this section, we give examples for degree $2$ divisors on $C_{B, z_o}$  which is a reduction of the degree $4$ divisor
given by the $(2,3, 6)$ quasi-torus decomposition of a $9$-cuspidal sextic $B$ considered in the previous section.
\begin{example}\label{ex:fib1_1}
\upshape
We consider the case $(a,b)=(1,1)$.
Then we have
\begin{align*}
	& g(t, x) = -\frac{1}{3} x^4 - 2 t x + 1, \quad
	h(t, x) = - \frac{2}{27} x^6 - \frac{2}{3} t x^3
		+ \frac{1}{3} x^2 - t^2 + 1, \\
	& f(t, x) = 4 g(t, x)^3 + 27 h(t, x)^2 \\
	&= -4 x^6 +4 t x^5 - x^4 + (4 t^3 - 36 t) x^3
	+ (30 t^2 + 18) x^2 - 24 t x + 27 t^4 - 54 t^2 + 31.
\end{align*}
Then $f(t,x)$ satisfies the condition described in
Section~\ref{subsec:double-cover}.
Let $\varphi_{B,z_o}\colon S_{B,z_o}\to\PP^1$ be
the genus 2 fibration as in Section~\ref{subsec:double-sextic}
such that the generic fiber $C_{B,z_o}$ is given by
$y^2=f(t,x)$.

Since $(y+\sqrt{27}h)(y-\sqrt{27}h)=4g^3$,
there exists a divisor $\ded_o$ such that
\[
	(y-\sqrt{27}h)
	= 3 \ded_o.
\]
Let $(u_0,v_0)$ be the Mumford representation of
the affine part of $\ded_o$.
Then we have
\[
	u_0 = -3 g, \quad
	v_0 = \sqrt{27} h \bmod u_0
	= - \frac{2}{\sqrt{3}} t x^3 + \frac{1}{\sqrt{3}} x^2
	- 3 \sqrt{3} \, t^2 + 3 \sqrt{3}.
\]
Applying Algorithm~\ref{alg:reduction}, we have
$\ded_o
\sim \ded_1 - \ded_{\infty}$, where
\begin{gather*}
	\ded_1 = \divisor[u_1, v_1], \quad
	u_1 = \frac{(t^2+3)x^2-4tx+1}{t^2+3}, \\
	v_1 = \frac{(26\sqrt{3}\,t^3-18\sqrt{3}\,t)x
	+9\sqrt{3}\,t^6+45\sqrt{3}\,t^4
	+20\sqrt{3}\,t^2-78\sqrt{3}}{3t^4+18t^2+27}.
\end{gather*}
Let $D$ be a divisor on $S_{B,z_o}$ with $D|_{C_S}=\ded_1$.
Then $D_{\mathrm{hor}}\cap(O^+\cup O^-)\neq\emptyset$
since $D_{\mathrm{hor}}$ and $O^+\cup O^-$ intersect
on the fibers at $t=\pm\sqrt{-3}$.

Note that the cusps of $B$ are $[1,0,0]$ and $[T_i,X_i,1]$, $i=1,2,\dotsc,8$,
where $T_1,T_2,\dotsc,T_8$ are the roots of
$27T^8+216T^6+756T^4-324T^2-676$, which has no multiple roots.
Any line through $z_o=[0,1,0]$ has
a equation of the form $\alpha T+\beta Z=0$.
Therefore any line through $z_o$ passes through at most one cusp of $B$
as in Section~\ref{subsec:proof}.
\end{example}
\begin{example}\label{ex:fib0_1}
\upshape
We consider the case $(a,b)=(0,1)$.
Then we have
\begin{gather*}
	g(t, x) = -\frac{1}{3} x^4 - 2 t x, \quad
	h(t, x) = - \frac{2}{27} x^6 - \frac{2}{3} t x^3
		- t^2 + 1, \\
	f(t, x) = 4 g(t, x)^3 + 27 h(t, x)^2
	= -4 x^6 + 4 t^3 x^3 - 36 t x^3 + 27 t^4 - 54 t^2 + 27.
\end{gather*}
Then $f(t,x)$ satisfies the condition described in
Section~\ref{subsec:double-cover}.
Let $\varphi_{B,z_o}\colon S_{B,z_o}\to\PP^1$ be
the genus 2 fibration as in Section~\ref{subsec:double-sextic}
such that the generic fiber $C_{B,z_o}$ is given by
$y^2=f(t,x)$.

Since $(y+\sqrt{27}h)(y-\sqrt{27}h)=4g^3$,
there exists a divisor $\ded_o$ such that
\[
	(y-\sqrt{27}h)
	= 3 \ded_o.
\]
Let $(u_0,v_0)$ be the Mumford representation of
the affine part of $\ded_o$.

Then we have
\[
	u_0 = -3 g, \quad
	v_0 = \sqrt{27} h \bmod u_0
	= -\frac{2}{\sqrt{3}} t x^3 - 3 \sqrt{3} \, t^2 + 3 \sqrt{3}.
\]
Applying Algorithm~\ref{alg:reduction}, we have
$\ded_o
\sim \ded_1 - \ded_{\infty}$, where
\[
	\ded_1 = \divisor[u_1, v_1], \quad
	u_1 = x^2, \quad
	v_1 = 3 \sqrt{3} (t^2 - 1).
\]
In other words, we have
\[
	\ded_1 = 2 (0, 3 \sqrt{3} (t^2 - 1)) - \ded_{\infty}. 
\]
Let $D$ be a divisor on $S_{B,z_o}$ with $D|_{C_S}=\ded_1$.
Then we have $D_{\mathrm{hor}}\cap(O^+\cup O^-)=\emptyset$.
Therefore the genus 2 fibration $\varphi_{B,z_o}$ does not satisfy
Theorem~\ref{thm:main}.

In fact, the cusps of $B$ are $[1,0,-1]$, $[1,0,0]$, $[1,0,1]$,
and $3$ points on each of the lines $T\pm\sqrt{3}\,Z=0$.
Since the lines $T\pm\sqrt{3}\,Z=0$ pass through $z_o=[0,1,0]$,
the point $z_o$ does not satisfy
the condition that any line through $z_0$ passes through
at most one cusp of $B$ assumed in Section~\ref{subsec:proof}.
\end{example}

\noindent Hiro-o TOKUNAGA, Yukihiro UCHIDA \\
Department of Mathematical Sciences, 
Tokyo Metropolitan University,   Hachiohji 192-0397 JAPAN , 
{\tt tokunaga@tmu.ac.jp, yuchida@tmu.ac.jp}


\begin{thebibliography}{99}
   
%
  \bibitem{bannai-tokunaga15} S.~Bannai and H.~Tokunaga: \emph{Geometry of bisections of elliptic surfaces and Zariski N-plets for conic arrangements},  Geom. Dedicata {\bf 178}(2015), 219 - 237.
  
   
 \bibitem{bpv} W.~Barth, K.~Hulek, C.A.M.~Peters and A. Van de Ven: Compact complex surfaces, 
 Ergebnisse der Mathematik und ihrer Grenzgebiete {\bf 4} 2nd Enlarged Edition, Springer-Verlag (2004).
 


\bibitem{cantor} D.~G.~Cantor: \emph{Computing in the Jacobian of a hyperelliptic curve}, Math. Comp. {\bf 48} (1987) 95--101.


%
\bibitem{cogo-lib}  J. -I.~Cogolludo-Agustin and A.~Libgober: \emph{Mordell-Weil groups of elliptic threefolds and the Alexander module of plane curves}, 
J. Reine Angew. Math. {\bf 697}(2014), 15-55.


%
\bibitem{fischer} G.~Fischer: {Plane Algebraic Curves}, Student Mathematical Library 15 (2001), American Math. Soc.
 
 %
 \bibitem{ghmm} S.~D.~Galbraith, M.~Harrison and D.~J.~Mireles Morales:
 \emph{Efficient hyperelliptic arithmetic using balanced representation for divisors. Algorithmic number theory}, 
 Lecture Notes in Comput. Sci., {\bf 5011}, Springer, (2008) 342 - 356.  
 
 %
 
 \bibitem{horikawa} E.~Horikawa: \emph{ On deformation of quintic surfaces},
\rm Invent. Math. {\bf 31} (1975), \rm $43 - 85$.

%


\bibitem{grant} D.~Grant: \emph{On an analogue of the Lutz-Nagell Theorem for hyperelliptic curves}, J. Number Theory {\bf 133} (2013),
963-969.

\bibitem{kloosterman2013} R.~Kloosterman: \emph{Cuspidal plane curves, syzygies and a bound on the MW-rank},  J. Algebra {\bf 375} (2013), 216-234

\bibitem{kloosterman2014} R.~Kloosterman: \emph{Mordell-Weil lattices and toric decompositions of plane curves},  Math. Ann. {\bf 367} (2017), 755-783.




\bibitem{lib} A.~Libgober:  \emph{On Mordell-Weil groups of isotrivial abelian varieties over function fields},  Math. Ann. {\bf 357} (2013), 605-629.

%
\bibitem{mwz} A.~J.~Menezes, Y.-H.~Wu and R.~J.~Zuccherato: \emph{An Elementary Introduction to Hyperelliptic 
Curves}, Appendix to N.~Koblitz, Algebraic Aspects of Cryptography, Algorithms and Computation in Math. {\bf 3}, Springer.
%

\bibitem{mumford} D.~Mumford: Tata Lectures on Theta II, Progress in Math. {\bf 43}, Birkh\"auser, Boston (1984).

%

\bibitem{oguiso-shioda}  K.~Oguiso and T.~Shioda: \emph{The Mordell-Weil lattice of Rational
Elliptic surface}, Comment. Math. Univ. St. Pauli \textbf{40}(1991), 83-99.

\bibitem{oka} M.~Oka:  \emph{Geometry of reduced sextics of torus type}, Tokyo J. Math. 26 (2003) 301-327.

\bibitem{oka-tai} M.~Oka and Pho Duc Tai: \emph{Classification of sextics of torus type}, Tokyo J. Math. {\bf 25} (2002) 399-433.


%
\bibitem{silverman} J.~H.~Silverman: {\rm The Arithmetic of Elliptic Curves}, Graduate Texts in Math. {\bf 106}, Springer (1986).
%



\bibitem{tokunaga94} H.~Tokunaga:  \emph{On dihedral Galois coverings}, \rm Canadian J. of
Math. {\bf 46} \rm (1994),1299 - 1317.


%
%
\bibitem{tokunaga99} H.~Tokunaga: \emph{Irreducible plane curves with the Albanese dimension 2}, Proc. AMS. {\bf 127}
(1999) 1935 - 1940.

\bibitem{tokunaga99-2} H.~Tokunaga: \emph{$(2, 3)$ torus sextic curves and the Albanese images of $6$-fold cyclic multiple planes},
 Kodai Math. J.{\bf 22} (1999), 222-242.
 
\bibitem{tokunaga02} H.~Tokunaga:  \emph{Galois covers for ${\mathfrak {S}}_4$ and ${\mathfrak {A}}_4$ and 
their applications}, Osaka Math. J., \textbf{39}(2002), 621-645. 
%
%

\end{thebibliography}
  \end{document}